\documentclass[11pt,a4paper]{article}
\usepackage{makeidx}
\usepackage{graphicx}
\usepackage{latexsym}
\usepackage{amsmath}
 \usepackage{amsfonts}
\usepackage{verbatim}
\usepackage{nameref}
\usepackage{hyperref}
\usepackage{sectsty}
\usepackage{enumerate}
\usepackage{subfigure}
\usepackage{epsf}
\usepackage{amssymb}

\usepackage{amsmath}
 \usepackage[dvips]{epsfig}
 \usepackage{pst-plot}
 \usepackage{psfrag}
\usepackage{a4wide}
\newcommand{\bX}{{\bf x}}
\newcommand{\WW}{{W}}
\newcommand{\hatW}{\hat{W}}
\newcommand{\barW}{V}

\newcommand{\alp}{\alpha}
\newcommand{\VV}{V}
\newcommand{\PP}{G}

\newcommand{\Oo}{\Omega}
\newcommand{\dO}{{\rm d} \Omega}
\newcommand{\dG}{{\rm d} \Gamma}
\newcommand{\dS}{{\rm d}S}
\newcommand{\dB}{{\rm d} \calB}

\newcommand{\Gu}{{\Gamma_{\chi}}}
\newcommand{\Gt}{{\Gamma_t}}
\newcommand{\St}{{S_t}}
\newcommand{\calB}{{\cal B}}
\newcommand{\uu}{{u}}

\newcommand{\be}{{\bf e}}
\newcommand{\cc}{{c}}
\newcommand{\xx}{{x}}
\newcommand{\la}{\label}
\newcommand{\half}{\frac{1}{2}}
\newcommand{\bt}{{\bf t}}

\newcommand{\bH}{{\bf H}}
\newcommand{\bE}{{\bf E}}
\newcommand{\bG}{{\bf G}}
\newcommand{\bI}{{\bf I}}
\newcommand{\bC}{{\bf C}}
\newcommand{\bQ}{{\bf Q}}
\newcommand{\bR}{{\bf R}}

\newcommand{\bF}{{\bf F}}
\newcommand{\bT}{{\bf T}}
\newcommand{\bn}{{\bf n}}
\newcommand{\bN}{{\bf N}}
\newcommand{\bx}{{\bf x}}

\newcommand{\sta}{{\rm sta}}
\newcommand{\vsig}{\zeta}

\newcommand{\bxi}{\mbox{\boldmath$\xi$}}
 
\newcommand{\bchi}{\mbox{\boldmath$\chi$}}

\newcommand{\barvsig}{{\bar{\varsigma}}}

\newcommand{\bsig}{\mbox{\boldmath$\sigma$}}
\newcommand{\bgamma}{\mbox{\boldmath$\gamma$}}
\newcommand{\btau}{{\mbox{\boldmath$\tau$}}}

\newcommand{\bbeta}{\mbox{\boldmath$\eta$}}

\newcommand{\calP}{{\cal P}}

\newcommand{\calE}{{\cal{E}}}
\newcommand{\calF}{{\cal{F}}}
\newcommand{\calS}{{\cal{S}}}
\newcommand{\barcalS}{\bar{\cal{S}}}
\newcommand{\calU}{{\cal{U}}}
\newcommand{\calT}{{\cal{T}}}
\newcommand{\calX}{{\cal{X}}}
\newcommand{\calG}{{\cal{G}}}

\newcommand{\real}{\mathbb{R}}
\newcommand{\eb}{\begin{equation}}
 \newtheorem{thm}{Theorem}
 \newtheorem{Lemma}{Lemma}
 \newtheorem{Corollary}{Corollary}
  \newtheorem{definition}{Definition}
\newcommand{\ee}{\end{equation}}
\newcommand{\calW}{{\cal W}}

\newcommand{\calR}{{\cal{R}}}

\newcommand{\SO}{\mbox{SO}}

\newcommand{\Lam}{{\Lambda}}
  \newcommand{\bS}{{\bf S}}
    \newcommand{\barbS}{\bar{\bf S}}
  \newcommand{\barbchi}{\bar{\bchi}}

 \newtheorem{rem}{Remark}

\def \tr{\mbox{tr\hskip 1pt}}

 \vspace{-1.5cm}
\title{\bf {\Large   Analytic Solutions to Large Deformation Problems Governed by
 Generalized Neo-Hookean Model  }}
 \author{David Yang  Gao    \\[0.2cm]
 \small   Federation University Australia, Mt Helen, VIC 3353, Australia }
\date{}
\begin{document}
\maketitle

\begin{abstract}
 This paper addresses some fundamental issues in nonconvex  analysis.
 By using pure complementary energy principle proposed by the author,
 a class of fully nonlinear partial differential equations   in nonlinear elasticity is able to converted a
 unified algebraic equation,
   a complete set of analytical solutions are  obtained for 3-D  finite deformation problems governed by
 generalized neo-Hookean model.
   Both global and local extremal solutions to the nonconvex variational problem  are identified by a triality theory.
 Connection between challenges in nonlinear  analysis and NP-hard problems in computational science is revealed.
Results show  that Legendre-Hadamard condition can only guarantee ellipticity for generalized convex problems.
For nonconvex systems, the ellipticity depends  not only on the stored energy, but also  on the external force field.
Uniqueness is proved based on a generalized  quasiconvexity and a generalized  ellipticity condition.
Application is illustrated for nonconvex logarithm stored energy.

\end{abstract}

 \section{ Nonconvex Variational Problem and Challenges}
Minimum total potential energy principle in nonlinear elasticity
has always presented fundamental challenging problems not only in
 continuum mechanics, but also in nonlinear analysis and computational sciences.
This paper intends to solve, under certain conditions,  the following   minimum potential variational problem  ($(\calP)$ for short):
  \eb
  (\calP): \;\; \min \left\{ \Pi(\bchi) = \int_{\calB} \WW(\nabla \bchi)  \dB - \int_{S_t} \bchi \cdot \bt \dS  | \;\; \bchi \in \calX_c \right\} , \label{eq-p}
  \ee
where the unknown deformation
  $\bchi (\bx) = \{ \chi_i (x_j) \} \in \calX_a$ is a vector-valued mapping
$\calB \subset \real^3 \rightarrow  \omega \subset  \real^3$
from a given material particle $\bx = \{ x_i\} \in \calB $ in the
undeformed body  to  a position vector  in the deformed configuration $ \omega$.
The body is fixed on the boundary $S_x \subset \partial \calB$, while on the remaining
boundary $S_t = S_x \cap \partial \calB$, the body is subjected to a given surface traction $\bt(\bx)$.
 In this paper, we let $\calX_a$ as a {\em geometrically admissible space} defined by
  \eb
  \calX_a = \{ \bchi \in \calW^{1,1} (\calB; \real^3) |  \;\;\bchi(\bx) = 0 \;\; \forall \bx \in S_x  \}
  \ee
where   $\calW^{1,1} $ is the standard notation for Sobolev space, i.e. a  function space in which both
$\bchi$ and its weak derivative  $\nabla \bchi$ have a finite $L^1(\calB)$ norm.
For homogeneous hyperelastic body, the  strain energy $\WW(\bF)$   is assumed to be $C^1$ on its domain
  $\calF_c \subset  \real^{3 \times 3} $, in which certain necessary  {\em constitutive constraints} are included, such as
\eb\label{const}
\det \bF > 0 , \;\;   \WW(\bF) \ge 0 \;\; \forall \bF \in \calF_c, \;\; \WW(\bF) \rightarrow \infty  \mbox{ as } \| \bF \| \rightarrow \infty.
\ee
Thus,   the {\em kinetically admissible space} in $(\calP)$ is simply defined by
\eb
\calX_c = \{ \bchi \in \calX_a | \; \nabla \bchi \in \calF_c \}
\ee
which is  essentially  nonconvex due to  nonlinear   constraints such as  $\det (\nabla \bchi)  > 0 $.
 Also, the stored energy $\WW(\bF)$ is  in general nonconvex in order to model real-world problems such as
 post-buckling and phase transitions, etc.
   Therefore,   the nonconvex variational problem $(\calP)$ has usually multiple local optimal solutions.

Let $\calX_b \subset \calX_c$ be a subspace with two additional conditions
\eb\label{cw}
\calX_b = \{ \bchi \in \calX_c | \; \; \bchi \in C^2(\calB; \real^3), \;\;  \WW(\bF(\bchi) )  \in C^2(\calF_c; \real) \}.
\ee
If $\partial \calB$ is sufficiently regular, the criticality condition $\delta \Pi(\bchi; \delta \bchi) = 0 \;\; \forall \delta \bchi \in \calX_b$ leads to  a nonlinear boundary-value problem
\eb
(BVP): \;\; \; \left\{
\begin{array}{l}
 - \nabla \cdot \bsig (\nabla \bchi) = 0 \;\; \mbox{ in } \calB, \\
\bN \cdot \bsig (\nabla \bchi) = {\bf t} \;\; \mbox{ on } S_t  ,\;\;   \bchi  = 0 \;\; \mbox{ on } S_x
 \end{array}\right.
 \label{eq-ebp}
 \ee
 where,  $\bN \in \real^3$ is a unit vector normal to $\partial \calB$,
  and $\bsig(\bF)$  is the first Piola-Kirchhoff   stress (force per unit undeformed
area),  defined by
\eb
\bsig   = \nabla \WW(\bF) , \;\; \mbox{ or } \;  \sigma_{ij}  = \frac{\partial \WW(\bF)}{\partial F_{ij} }, \;\;
i,j = 1,2,3.
\ee

\begin{rem}[Nonconvexity, Multi-Solutions, and NP-Hard Problems]$\;$\newline
{\em
The stored energy $\WW(\bF)$ in nonlinear elasticity is generally nonconvex. It turns out that the fully nonlinear $(BVP)$
could have  multiple solutions $\{\bchi_k(\bx) \}   \in \calX_c \subset \real^\infty$ at each material point $\bx\in \calB_s \subset \calB$.
As long as the continuous domain  $\calB_s \neq \emptyset$, this solution set $\{ \bchi_k(\bx) \} (k=1, \dots , K)$ can form  infinitely many ($K^\infty$) solutions
to $(BVP)$ even $\calB \subset \real$. It is impossible to use traditional convexity and ellipticity conditions to identify global minimizer
among all these local solutions.
 Gao and Ogden discovered in \cite{gao-ogden-qjmam} that
for certain given external force field, both global and local extremum solutions are nonsmooth and can't be obtained
by  Newton-type numerical methods. Therefore, Problem $(\calP)$ is much more difficult than $(BVP)$.
 In computational mechanics, any direct numerical method  for solving   $(\calP)$ will lead  to
a nonconvex minimization problem in $\real^n$, which could have $K^n$ local solutions.
Due to the lack of global optimality condition, it is fundamentally difficult to solve nonconvex minimization problems by traditional methods within polynomial time. Therefore,  in  computational sciences most nonconvex minimization problems are  considered to be NP-hard (Non-deterministic Polynomial-time hard)  \cite{gao-bridge}.

Direct methods for solving nonconvex  variational problems   in finite elasticity have been
studies extensively during the last fifty years and  many generalized
convexities, such as poly-, quasi- and rank-one convexities, have been proposed.
For a  given function  $W:\calF_c \rightarrow \real$, the
following statements  are well-known (see  \cite{sch-neff})\footnote{It was proved recently that rank-one convexity  also implies polyconvexity for
isotropic, objective and isochoric elastic
energies in the two-dimensional case \cite{neff15}.}:
\[
 \mbox{    convex $\Rightarrow   \mbox{ poly-convex }  \Rightarrow
  \mbox{ quasi-convex }   \Rightarrow   \mbox{ rank-one convex}$.}
  \]

Although the generalized convexities have been well-studied  for general function $\WW(\bF)$ on matrix space
$\real^{m\times n}$, these mathematical concepts  provide only necessary conditions for local minimal solutions, and
can't be   applied  to general (nonconvex) finite deformation problems.
In reality, the stored energy $\WW(\bF)$ must be nonconvex in order to model real-world phenomena.
Strictly speaking,  due to certain necessary  constitutive constraints such as  $\det \bF > 0$  and objectivity
condition etc,
even the domain $\calF_c$ is not convex, therefore, it is not appropriate to discuss convexity of the stored energy $\WW(\bF)$ in general nonlinear elasticity.
  How to identify global optimal solution has been a fundamental challenging problem in nonconvex analysis and computational science.
\hfill $\blacksquare$ }
\end{rem}

\begin{rem}[Canonical Duality,  Gap Function, and Global Extremality]$\;$\\
{\em
The objectivity is a necessary constraint   for any hyper-elastic model.
 A real-valued function  $\WW:\calF_c \rightarrow \real$  is  objective  iff
 there exists a function $\VV(\bC)$ such that
$
W(\bF)  = \VV(\bF^T \bF) \;\; \forall \bF \in \calF_c
 $ (see \cite{ciarlet}).
By the fact that the right Cauchy-Green tensor $\bC$ is an objective measure on a convex domain
$\calE_a = \{ \bC \in \real^{3\times 3} | \;\; \bC = \bC^T, \;\; \bC \succ 0 \}$, it is possible and natural
to discuss the convexity of
$\VV(\bC)$. A real-valued function $\VV:\calE_a \rightarrow \real$ is called {\em  canonical} if the duality relation
 $\bxi^* = \nabla \VV(\bxi) : \calE_a \rightarrow \calE^*_a$ is one-to-one and onto \cite{gao-dual00}.
 The  canonical duality is necessary for modeling natural phenomena, which  lays a foundation for the canonical duality theory \cite{gao-dual00}. This theory was developed from
  Gao and Strang's original work  in 1989 \cite{gao-strang89a} for general nonconvex/nonsmooth  variational problems in finite deformation theory.
  The key idea of this theory is assuming the existence of a  geometrically admissible (objective) measure
  $\bxi=\Lambda(\bF)$ and a canonical function $\VV(\bxi)$ such that the following {\em  canonical transformation}  holds
 \eb
 \bxi = \Lambda(\bF):\calF_a \rightarrow \calE_a   \;\;  
 \Rightarrow \; \; \WW(\bF) = \VV(\Lambda(\bF)). \label{eq-lam}
 \ee
  Gao and Strang discovered that  the  directional  derivative $\Lam_t (\bF)= \delta\Lam(\bF) $  is adjoined with
the  equilibrium operator,  while its complementary operator $\Lam_c(\bF) = \Lambda(\bF) - \Lam_t(\bF) \bF $   leads to
 a so-called {\em    complementary gap function}, which recovers  duality gaps in traditional duality theories  and  provides a sufficient  condition for identifying both global and local extremal solutions
 \cite{gao-dual00,gao-bridge}. $\;$\hfill $\blacksquare$
 }\end{rem}

The canonical duality  theory has been applied for solving a large class of nonconvex, nonsmooth, discrete problems
  in multidisciplinary fields of nonlinear  analysis,   nonconvex  mechanics, global optimization, and computational sciences, etc.
  A comprehensive review is given recently in \cite{gao-bridge}.
 The main goal of this paper is to show author's  recent  analytical  solutions \cite{gao-cmt15}
  for general anti-plane shear problems   can be easily generalized for solving finite deformation problems governed by generalized neo-Hookean materials.
  Some insightful  results  are obtained on generalized convexity and ellipticity   in  nonlinear analysis.

\section{Complete Solutions to Generalized Neo-Hookean Material}
By the fact that the  right Cauchy-Green strain $\bC = \bF^T \bF$ is  an objective tensor,  its three
  principal invariants
\begin{equation}
I_1(\bC)=\tr\mathbf{C},\quad  I_2 (\bC)= \half [ (\tr \bC)^2 - \tr(\bC^2) ]
 ,\quad I_3 (\bC)=\det\mathbf{C}
\end{equation}
are also objective functions of $\bF$. Clearly, for  isochoric   deformations we have $ I_3 (\bC) = 1 $.
The elastic body is said to be  {\em  generalized neo-Hookean material} if  the stored energy depends only on  $I_1$,
  i.e. there exists a function $\barW(I_1) $
 such that $W(\bF) =  \barW(I_1(\bC(\bF)) ) $.
Since $I_1 = \tr(\bF^T\bF) > 0 \;\; \forall \bF \in \calF_c$, the   domain of $\barW(I_1)$
   is a   convex (positive)  cone
   \eb
   \calE_a =   \{ \xi \in L^p(\calB) \; | \; \xi(\bx)  > 0  \;\; \forall \bx \in \calB \},
   \ee
   it is possible to discuss  the  convexity of $\barW(I_1)$ on $\calE_a$.
  Furthermore, we assume that $\barW(I_1)$ is a $C^2(\calE_a )$ canonical function.
   Then the canonical transformation (\ref{eq-lam}) for the generalized neo-Hookean model is
\eb
 \xi = \Lam(\bF) = \tr(\bF^T \bF) : \calF_c \rightarrow \calE_a, \;\;
 \WW(\bF) = \barW(\xi(\bF)) .
\ee
For a  given external force $\bt(\bx) $ on $\St$, we introduce
  a {\em statically admissible space}
\eb
{\calT}_a = \left\{ { \bT } \in \calW^{1,1}(\calB ; \real^{3\times 3}) \; | \;\; \nabla \cdot \bT = 0 \;\; \mbox{ in } \calB, \;\; \bN \cdot \bT =  \bt \;  \;\; \mbox{ on } \St \right\}.
 \ee
Thus for any given $\bT \in \calT_a$, 
the primal problem $(\calP)$ for the generalized neo-Hookean material can be written in following canonical form
 \eb
 (\calP)_{\bT} : \;\; \min \left\{ \Pi_\bT (\nabla \bchi) = \int_{\calB} \PP(\nabla \bchi)\; \dB  \;\;   | \;\;
 \forall \bchi \in \calX_c \right\},
 \ee
where $\calX_c = \{ \bchi \in \calX_a| \;  \Lambda(\nabla \bchi) \in \calE_a \}$ and the integrand
$\PP:\calF_a \rightarrow \real$ is defined by
\eb
 \PP(\bF) = \barW(\Lambda(\bF)) - \tr(\bF^T \bT ).
 \ee
By the fact that $\det \bF > 0$ is not a variational constraint and
the certain constitutive constraints, such as coercivity and   objectivity,
have been naturally relaxed by the canonical transformation, the domain
of $\PP(\bF)$  is simply $\calF_a =  \real^{3\times 3}$.

 Let $\SO(3) = \{ \bR \in \real^{3\times 3} | \;\; \bR^T  = \bR^{-1}, \;\; \det \bR = 1 \}$ and
\eb
\calR = \{ \bR(\bx)  \in L^1 [ \calB , \real^{3\times 3}] | \;\; \bR(\bx) \in \SO(3) \;\; \forall \bx \in \calB \}.
\ee

\begin{thm}\label{thm1}
For any given $\bT \in \calT_a$, if $\barbchi \in \calX_c$ is a stationary solution to  $(\calP)_\bT$, then it is also a stationary  solution to $(\calP)$.

 For any given rotation field $\bR(\bx) \in \calR$ such that $\bR^T \bT \in \calT_a$, then  $\Pi_{\bT}(\bF) = \Pi_{\bT}(\bR \bF)$.

 For any uniform rotation  $\bR \in \SO(3)$ such that $\bR^T \bT \in \calT_a$,  if
$\barbchi$ is a stationary solution to  $(\calP)$, then $\bR \barbchi$ is also a stationary  solution to $(\calP)$.
\end{thm}
{\bf Proof}. For any given $\bT \in \calT_a$, the stationary condition for the canonical variational problem $(\calP)_\bT$ leads to the following canonical
boundary value problem
\eb\label{bvpg}
(BVP)_{\bT} : \;\;
\left\{ \begin{array}{l}
\nabla \cdot (2 \vsig  \nabla \bchi)  =  \nabla \cdot \bT =   0 \;\; \mbox{ in } \calB, \\
\bN \cdot (2 \vsig  \nabla \bchi)  = \bN \cdot \bT =  \bt \; \mbox{ on } \St, \;\; \bchi = 0 \;\; \mbox{ on } S_x
\end{array} \right.
\ee
which is identical to $(BVP)$ since
\[
\bsig = \nabla \WW(\bF) =  \frac{\partial \VV(\xi)}{\partial \xi} \frac{\partial \xi}{\partial \bF}
   = 2 \vsig  \bF,  \;\;  \vsig  =  \nabla \VV(\xi)  . 
  \]

  By the objectivity of   $\xi =\Lambda(\bF) = \Lambda(\bR \bF)  \;\; \forall \bR(\bx)  \in \calR $ and the fact that
  \[
  \int_{\calB} \tr [(\bR \nabla \bchi)^T \bT ] \dB =      \int_{\calB} \tr [(\nabla \bchi)^T (\bR^T \bT) ] \dB  = \int_\St \bchi \cdot \bt \dS \;\; \forall \bR^T \bT \in \calT_a,
  \]
  we have $\Pi_{\bT}(\bF )  = \Pi_\bT(\bR \bF) \;\; \forall \bR (\bx) \in \calR$.
  Particularly,   for any uniform  $\bR \in  \SO(3) $ such that $\bR^T \bT \in \calT_a$, we have $\Pi(\bchi) = \Pi_{\bT}(\bR \bF(\bchi))$. \hfill $\Box$\\

 Theorem \ref{thm1} is important for understanding the canonical duality theory.

By the canonical assumption on $\barW(\xi)$, the duality relation   $\vsig = \nabla \barW(\xi) : \calE_a \rightarrow \calE^*_a$ is invertible.
The complementary energy  can be defined uniquely  by the
Legendre transformation
\eb
\barW^*(\vsig) = \{\xi \vsig - \barW(\xi) | \; \vsig = \nabla \barW(\xi) \}.
\ee
 Clearly, the function   $\barW:\calE_a \rightarrow \real$ is   canonical if and only if the following
 {\em canonical duality relations} hold on $\calE_a \times \calE^*_a$
\eb\label{cdr}
\vsig = \nabla \barW(\xi) \;\; \Leftrightarrow \;\; \xi = \nabla \barW^*(\vsig)  \;\; \Leftrightarrow
\;\;
\barW(\xi) + \barW^*(\vsig) = \xi \vsig.
\ee
Using $\barW(\xi) = \xi\vsig - \barW^*(\vsig)$, the nonconvex function
$\PP(\bF)$ can be written as the standard Gao and Strang    total complementary function
 $\Xi: \calX_a \times \calE^*_a \rightarrow \real $
\eb
\Xi(\bchi, \vsig) =  \int_{\calB} \left[ \Lam(\nabla \bchi) \vsig - \barW^*(\vsig) - \tr( (\nabla \bchi)^T \bT ) \right] \dB.
\ee
Let $\calS_a \subset \calE^*_a$ be a canonical dual feasible space defined by
 \eb
 \calS_a = \{ \vsig \in \calE^*_a | \; \vsig^{-1} \tau^2 \in L^1(\calB) \}.
 \ee
 Then for a given $\vsig \in \calS_a$, the canonical dual function
  can be obtained by the {\em canonical dual transformation}:
\eb\label{eq-pdxi}
\Pi^d(\vsig) = \sta \{ \Xi(\bchi, \vsig) | \;\; \bchi \in \calX_a  \}
=   \int_{\calB} \PP^d(\vsig) \dB,
  \ee
  where the notation $\sta \{ \Xi(\bchi, \vsig)| \;\bchi \in \calX_a\}$ stands for finding (partial)
   stationary point $\bchi \in \calX_a$ of $\Xi(\bchi, \vsig)$ for a given $\vsig \in \calS_a$, and
  \eb
   \PP^d(\vsig) =     - \barW^*(\vsig) - \frac{1}{4} \vsig^{-1} \tau^2   ,  \;\;  \;\; \tau^2 =   \tr(\bT^T  \bT ) .
   \ee
 Thus,  the pure complementary energy  principle, first proposed   in 1998 \cite{gao-ima98},
 leads to the following canonical dual variational problem
  \eb
 (\calP^d) : \;\;\;\; \sta \left\{  \Pi^d(\vsig) =   \int_{\calB} \PP^d(\vsig) \dB \; | \;\; \vsig \in \calS_a \right\} .
  \ee
 Since the canonical dual variable $\vsig$ is a scalar-valued function, the criticality condition
 for this variational problem  leads to  a
so-called   {\em canonical dual algebraic equation} (see \cite{gao-dual00}):
  \eb\label{cda}
  4 \vsig^2  \nabla \barW^*(\vsig)  =  \tau^2(\bx) \;\; \;\; \forall    \bx \in  \calB.
  \ee
Note that  $\nabla \barW^*(\vsig):\calE^*_a \rightarrow \calE_a$ is also one-to-one and onto,
 this equation has at least one solution for any given $\tau^2 = \tr(\bT^T\bT) \ge 0 $ and $\vsig = 0 $ only if $\tau = 0$. Therefore, although there is an inverse term $\zeta^{-1}$ in  $\PP^d(\zeta)$, this
  canonical dual function is well-defined on  $\calS_a$.
Due to the nonlinearity, the solution to (\ref{cda}) may not be unique \cite{gao-dual00,gao-cmt15,gao-ogden-qjmam}.
By the pure complementary energy principle proposed by  Gao in 1999
(see \cite{gao-dual00}), we have
\begin{thm}[Complementary-Dual  Principle]
For any given  $\bT \in \calT_a$, the following statements are equivalent:

1)  $(\barbchi, \barvsig)$ is a stationary point of $\Xi(\bchi, \vsig)$;

2) $\barbchi$ is a stationary solution to $(\calP)$;

3) $\barvsig$ is a stationary solution to $(\calP^d)$.

Moreover, we have
\eb\label{eq-pdp}
\Pi(\barbchi ) = \Xi(\barbchi, \barvsig) =  \Pi^d( \barvsig )
\ee
 \end{thm}
 {\bf Proof}. For any given  $\bT \in \calT_a$, the stationary condition of $\Xi(\bchi, \vsig)$ leads to the canonical equilibrium equations
 \begin{eqnarray}
 & & \Lambda(\bF(\barbchi)) = \nabla \VV^*(\barvsig) , \label{eq-cdi}\\
 & & 2 \barvsig \bF(\barbchi) = \bT \in \calT_a  \label{eq-cec}
 \end{eqnarray}

  By the canonical duality,   (\ref{eq-cdi}) is  equivalent to $\barvsig = \nabla \VV(\xi)$ with $\xi = \Lam(\nabla \barbchi)$.
Thus, $\barbchi$ must be a stationary solution to $(\calP)_{\bT}$ and also a stationary solution to $(\calP)$ due to Theorem \ref{thm1}.

By solving (\ref{eq-cec}) we have    $\bF(\barbchi) = \frac{1}{2  \barvsig} \bT$.
Substituting this into (\ref{eq-cdi}) leads to the canonical dual equation (\ref{cda}).
Thus, $\barvsig$ is a stationary solution to $(\calP^d)$.

The equivalence and
the equation  (\ref{eq-pdp}) can be proved  by
\[
\sta \{ \Pi_{\bT}(\nabla \bchi) | \; \bchi \in \calX_c\} = \sta \{  \Xi(\bchi, \vsig)
| \;\; (\bchi,   \vsig)  \in  \calX_a \times  \calE^*_a\} =   \sta \{ \Pi^d(\vsig) | \;\; \vsig \in \calS_a\}
\]
and Theorem \ref{thm1}.   \hfill $\Box$

 \begin{thm}[Pure Complementary Energy Principle]
For any given nontrivial $\bt \neq 0$ and $\bchi \in \calX_a$ such that  $\bT \in \calT_a \neq \emptyset$,  (\ref{cda}) has at least one solution    $\vsig_k \neq 0$,
 the deformation gradient  defined by $\bF_k = \nabla \bchi_k = \vsig_k^{-1}  \bT$
   is a critical point of $\Pi(\bchi)$ and
 $\Pi(\bchi_k) = \Pi^d( \vsig_k) $.

Moreover, if  $\nabla \times ( \vsig_k^{-1}  \bT) = 0$, then the deformation vector  defined by
   \eb\label{eq-solu}
   \bchi_k (\bx) = \half \int_{\bx_0}^\bx  \vsig_k^{-1}  \bT\cdot \mbox{d}\bx
   \ee
   along any path from $\bx_0 \in  S_x$ to $\bx \in \calB$
   is  a solution to $(BVP)_\bT$
   in the sense that it satisfies both equilibrium equation and boundary conditions in (\ref{bvpg}).
   \end{thm}
{\bf Proof}.
By the canonical duality relations in (\ref{cdr}) we know that  $\xi_k = \nabla \VV^*(\zeta_k) > 0$.
Thus, for a given  nontrivial $\bt(\bx)$, there exists a nontrivial  $\tau^2(\bx) = \tr(\bT^T \bT)$ in $\calB$ such that
the canonical dual algebraic equation (\ref{cda}) have at least one nontrivial solution $\zeta_k(\bx)$ in $\calB$.

 Since the critical point $\zeta_k $ is a solution to (\ref{cda}), we have
 \eb
 \xi_k = \tr (\bF^T_k \bF_k) = \frac{1}{4} \vsig_k^{-2} \tr(\bT^T \bT) = \nabla \VV^*(\vsig_k) \;\; \Rightarrow
  \;\;\bF_k = \half \vsig_k^{-1} \bT
 \ee
 subjected to any given  rotation field  $\bR(\bx)  \in \calR$. By the fact that the canonical dual solution $\vsig_k$
 defined by (\ref{cda})  is independent of the rotation field, the  canonical duality leads to
 \[
\PP^d(\vsig_z) = \Xi(\bF_k, \zeta_k)  = \VV(\Lambda(\bF_k)) - \tr(\bF_k^T \bT)= \PP(\bF_k).
 \]
 This shows $\Pi(\bchi_k) = \Pi^d(\zeta_k)$.

 To prove  $\bchi_k$ defined by (\ref{eq-solu}) is a solution to $(BVP)_\bT$, we simply substitute $\nabla \bchi_k = \bF_k  = \half \zeta_k^{-1} \bT$
 into $(BVP)_\bT$ to have all  necessary equilibrium conditions satisfied. Therefore, $\bchi_k$ defined by (\ref{eq-solu})
 is a solution to $(BVP)_\bT$.  \hfill $\Box$\\

 This  pure complementary energy  principle
 shows that by the canonical dual transformation, the fully nonlinear partial differential equation
 in $(BVP)_\bT$  can be converted to an algebraic equation (\ref{cda}),
 which can be solved   to obtain a complete set of solutions (see \cite{gao-cmt15,gao-haj}).

  Since $\calS_a $ is nonconvex, in order to identify global and local optimal solutions, we need
 the following convex subsets
\eb
\calS_a^+ = \{ \vsig \in \calS_a| \; \zeta   > 0 \}, \;  \;\;
 \calS_a^- = \{ \vsig \in \calS_a| \;  \zeta  < 0
 \}.
\ee
Then by  the canonical duality-triality theory developed in \cite{gao-dual00} we have   the following theorem.

\begin{thm}\label{thm-1}
Suppose that $\VV:\calE_a \rightarrow \real$ is convex and  for a  given $\bT \in \calT_a$ such that $\{\vsig_k\}$ is a solution set to (\ref{cda}),
 $\bF_k =  \half  \vsig^{-1}_k \bT  $, and
 $\bchi_k$ is defined by (\ref{eq-solu}),
 we have the following statements.
\begin{verse}
1. If $\vsig_k \in \calS_a^+$,
then  $\nabla^2 \WW(\bF_k) \succ  0$ and   $\bchi_k$ is a global minimal solution to $(\calP) $.\\

2. If $\vsig_k \in \calS_a^-$ and $\nabla^2  \WW(\bF_k)\succ 0$, then $\bchi_k$ is a local minimal solution to  $(\calP) $.\\

3. If $\vsig_k\in \calS_a^-$ and $\nabla^2  \WW(\bF_k) \prec 0$, then $\bchi_k$ is a local maximal solution to  $(\calP) $.
\end{verse}

If $\{\vsig_k\} \subset  \barcalS^+_a $,  then $\{\bchi_k\}$ is a convex set.
The solution of $(\calP) $  is unique  if $\{\vsig_k\} \subset \calS^+_a$. 
\end{thm}
{\bf Proof}.
 By using  chain rule for $\WW(\bF ) = \barW(\xi(\bF))$ we have $\nabla \WW(\bF) = 2 \bF   [\nabla \barW(\xi)] = 2 \vsig  \bF $, and
 \eb \label{eq-hatw}
   \nabla^2 \WW(\bF ) = 2 \vsig  \bI\otimes \bI   + 4  h(\xi)  \bF   \otimes \bF  ,
 \ee
  where $\mathbf{I}$ is an identity tensor in $\real^{3\times 3}$,  $h(\xi) =   \nabla^2 \barW(\xi) \ge 0 $ due to the convexity of   $\barW$
 on $ {\calE_a}$.
 Therefore, $ \nabla^2 \WW(\bF_k )  \succ  0$ if $\vsig_k \in \calS^+_a$.

 To prove $\bchi_k$ is a global minimizer of $(\calP)$, we  follow  Gao and Strang's work in 1989 \cite{gao-strang89a}. By the
  convexity of $\barW(\xi)$ on its convex domain $\calE_a$, we have
\eb\label{eq-xic}
\barW(\xi) - \barW(\xi_k) \ge   (\xi - \xi_k)  \vsig_k \;\;  \;\; \forall \xi,\; \xi_k \in \calE_a,\;\; \vsig_k = \nabla \barW(\xi_k).
\ee
For any given variation $\delta \bchi$, we let $\bchi= \bchi_k + \delta \bchi$.
Then we have \cite{gao-strang89a}
\eb\label{lamtc}
\Lambda(\nabla\bchi) = \tr[(\nabla \bchi)^T    (\nabla \bchi)]  = \Lambda(\nabla\bchi_k) + \Lam_t (\nabla \bchi_k)   (\nabla \delta  \bchi)    - \Lambda_c(\nabla \delta \bchi),
\ee
where $\Lam_t(\bF) \delta \bF = 2 \tr[ \bF^T  (\delta \bF) ]$ and  $
\Lambda_c(\delta \bchi)  =  - \Lambda ( \delta  \bchi)$. Clearly, $\Lam(\bF) = \Lam_t(\bF) \bF + \Lam_c(\bF)$.
Then combining  the inequality (\ref{eq-xic}) and (\ref{lamtc}), we have
\begin{eqnarray}
\Pi(\bchi) - \Pi(\bchi_k) &\ge& \int_\calB 2 \vsig_k \tr[ (\nabla \bchi_k)^T   (\nabla \delta \bchi)]   \dB
 - \int_\St \delta \bchi \cdot \bt \dS  +  \int_\calB  \vsig_k \tr [(\nabla \bchi)^T  (\nabla \bchi) ] \dB\nonumber \\
 &=& \int_\calB    [2\vsig_k  (\nabla \bchi_k)  - \bT ] : (\nabla \delta \bchi)   \dB
 + G_{ap}(\delta\bchi, \vsig_k)  \;\; \forall  \bchi, \; \delta \bchi \in \calX_c
\end{eqnarray}
for any given $ \bT \in \calT_a$, where
\eb
G_{ap}(\bchi, \vsig) = \int_\calB -\Lambda_c(\nabla\bchi) \vsig \dB = \int_\calB  \vsig \tr [(\nabla \bchi)^T  (\nabla \bchi) ] \dB
\ee
is the {\em complementary gap function} introduced by Gao and Strang in \cite{gao-strang89a}.
If $\bchi_k$ is a critical point of $\Pi(\bchi)$, then we have
\[
\int_\calB   [2 (\nabla \bchi_k) \vsig_k - \bT ] : (\nabla \delta \bchi)   \dB = 0 \;\; \forall  \delta \bchi \in \calX_c, \;\;
\forall \bT \in \calT_a
\]
Thus, we have
$
\Pi(\bchi) - \Pi(\bchi_k) \ge G_{ap}(\delta \bchi, \vsig_k) \ge 0 \;\; \forall  \delta \bchi \in \calX_c  \;  \mbox{ if } \vsig_k \in \calS_a^+.
$
This shows that $\bchi_k$ is a global minimizer of $(\calP) $.

To prove the local extremality, we replace $\bF_k$ in  (\ref {eq-hatw}) by  $\bF_k = \half \vsig^{-1}_k \bT$ such that
 \eb
 {\bf G} (\zeta_k) =   \nabla^2 \WW(\bF_k ) = 2 \vsig_k \bI\otimes \bI   + \vsig_k^{-2}  h(\xi_k)  \bT   \otimes \bT  ,
 \ee
where $\xi_k = \nabla \barW^*(\vsig_k)$. Clearly, for a given $\bT \in \calT_a$ such that $\vsig_k \in \calS^-_a$,
 the Hessian $ \nabla^2 \WW(\bF_k)$ could be either positive or negative definite.
The total potential $\Pi(\bchi_k)$ is locally convex if the {\em Legendre condition}
$ \nabla^2 \WW(\nabla \bchi_k) \succeq 0$ holds, locally concave
 if   $ \nabla^2 \WW(\nabla \bchi_k) \prec 0$.
 Since $\bchi_k$ is  a global minimizer when $\vsig_k \in \calS^+_a$,
 therefore, for   $\vsig_k \in \calS^-_a$, the stationary solution $\bchi_k$ is  a local  minimizer  if
 $ \nabla^2 \WW(\nabla \bchi_k) \succ  0$ and, by the triality theory\cite{gao-dual00,gao-bridge},
 $\bchi_k$ is the biggest local maximizer if
 $ \nabla^2 \WW(\nabla \bchi_k) \prec 0$.

 If $\{\vsig_k \} \subset  \calS^+_a $, then  all    the solutions $\{ \bchi_k \} $ are global minimizers and
 form a convex set. Since $\Pi^d(\zeta)$ is strictly concave on the open convex set $\calS^+_a $,
the condition  $\{\vsig_k \} \subset \calS^+_a$ implies the unique solution of (\ref{cda}).
 In this case, Problems $(\calP)_{\bT}$ has  at most one  solution.
 \hfill $\Box$
\begin{thm}[Triality Theory]
For any given $\bT \in \calT_a\neq\emptyset$,  let  $\zeta_k$ be    a critical point of   $ (\calP^d) $,
the vector   $\bchi_k$ be  defined by (\ref{eq-solu}), and
 $ \calX_o \times \calS_o  \subset
 \calX_c \times \calS^-_a$   a neighborhood\footnote{The neighborhood $\calX_o$ of $\bchi_k$ in the canonical duality theory  means that $\bchi_k$ is the only one critical point of $\Pi(\bchi)$ on $\calX_o$ (see \cite{gao-dual00}).}
 of $ (\bchi_k, \zeta_k)$.

If $\zeta_k \in \calS^+_a$, then
 \eb
 \Pi(\bchi_k) = \min_{\bchi \in \calX_c} \Pi(\bchi) = \max_{\vsig \in \calS^+_a} \Pi^d (\vsig) = \Pi^d(\vsig_k).
 \ee

If $\zeta_k \in \calS^-_a$ and  $ {\bf G} (\zeta_k) \succ 0$, then
\eb
 \Pi(\bchi_k) = \min_{\bchi \in \calX_o} \Pi(\bchi) = \min_{\vsig \in \calS_o} \Pi^d (\vsig) = \Pi^d(\vsig_k).
 \ee

 If $\zeta_k \in \calS^-_a$  and $ {\bf G} (\zeta_k) \prec  0$,
 then
 \eb
 \Pi(\bchi_k) = \max_{\bchi \in \calX_o} \Pi(\bchi) = \max_{\vsig \in \calS_o} \Pi^d (\vsig) = \Pi^d(\vsig_k).
 \ee
 \end{thm}

This theorem   shows that for convex canonical function $\VV$, the triality theory can be used to identify
both global and local extremum solutions to  the  variational problem $(\calP) $ and
the nonconvex minimum variational problem $(\calP)_{\bT}$ is canonically equivalent to
the following concave maximization problem over an open  convex
set $\calS^+_a$, i.e.
 \eb
 (\calP^\sharp)_{\bT}: \;\;\;\; \max \left\{  \Pi^d(\vsig) =   \int_{\calB} \PP^d(\vsig) \dB \; | \;\; \vsig \in \calS^+_a \right\} ,
  \ee
  which is much easier to solve than  directly  for obtaining  global optimal solution of $(\calP) $.

  \section{Generalized Quasiconvexity, G-Ellipticity,  and Uniqueness}

Ellipticity is a classical concept originally from linear partial differential systems,
 where  the deformation is a scalar-valued function $\chi:\calB \rightarrow \real$ and  stored energy is a  quadratic function
$\WW(\bgamma) = \half \bgamma^T  {\bH}  \bgamma $  of $\bgamma = \nabla \chi \in \real^3$.
The linear operator
\[
L[\chi] = - \nabla \cdot  [\bH   (\nabla \chi)]  =  - [h_{ij} \chi_{,j}]_{,i}
\]
 is called elliptic if
 $\bH = \{ h_{ij} \} $ is positive definite.  In this case, the function  $\PP(\bgamma ) = \WW(\bgamma)  - \bgamma^T \btau$
 is convex  and its level set  is an ellipse for any given $\btau \in \real^3$.
This concept has been extended  to   nonlinear analysis.
The fully nonlinear  partial differential  equation
 in $(BVP)$ (\ref{eq-ebp}) is called elliptic   if
  the following Legendre-Hadamard (LH) condition holds
  \eb \label{eq-LH}
({\bf a} \otimes {\bf a} ) :  \nabla^2\WW(\bF) : ({\bbeta}  \otimes {\bbeta} ) \ge 0 \;\; \forall {\bf a}, \bbeta \in \real^3,
 \;\; \forall \bF \in \calF_a .
\ee
The $(BVP)$ is called strong elliptic if the inequality holds strictly. In this case, $(BVP)$ has at most one solution.
In vector space,   the  LH condition is equivalent to Legendre condition $\nabla^2 \WW(\bgamma) \succeq 0  \;\; \forall \bgamma \in \real^n$.

Clearly, the LH condition is only a sufficient condition  for local minimizer of the variational problem $(\calP)$.
 In order to identify ellipticity, one must to check LH condition
for all local solutions, which is impossible for general fully nonlinear problems.
 Also, the traditional  ellipticity  definition  depends only on the stored energy $\WW(\bF)$
regardless of the linear term in $\PP(\bF) = \WW(\bF) - \tr(\bF^T \bT)$.
This definition works only for convex systems since the linear term $\tr(\bF^T \bT)$  can't change the convexity of $\PP(\bF)$.
But this is not true for  nonconvex  systems.
To see this, let us consider the St. Venant-Kirchhoff material
\eb
W(\bF) = \half \bE : {\bf H} : \bE, \;\; \bE = \half [(\bF)^T(\bF) - \bI ],
\ee
where $\bI$ is  a unit tensor in $\real^{3 \times 3}$. Clearly, this function is not even rank-one convex.
A special case of this model in $\real^n$  is the well-known double-well potential
$\WW(\bgamma) = \half ( \half |\bgamma|^2- 1)^2$.
If we let $\xi = \Lambda(\bgamma) = \half   |\bgamma|^2 - 1$ be  an objective measure, we have the
canonical function $\barW(\xi) = \half \xi^2 $. In this case, the canonical dual algebraic equation
(\ref{cda}) is  a cubic  equation (see \cite{gao-dual00})
$2 \zeta^2 (\zeta +1) = \tau^2$, which has at most three real solutions $\{\zeta_k(\bx) \}$ at each $\bx \in \calB$
satisfying $\zeta_1 \ge 0 \ge \zeta_2\ge \zeta_3$.
It was proved in \cite{gao-dual00} (Theorem 3.4.4, page 133)
that for a given force $\bt(\bx)$,
if  $\tau^2(\bx) > 8/27 \;\; \forall \bx \in \calB \subset \real$, then $(BVP)_{\bT}$ has only one solution on $\calB$.
If $\tau^2(\bx) < 8 /27 \;\; \forall \bx \in \calB_s  \subset \calB$, then
 $(BVP)_{\bT}$
has   three solutions $\{\chi_k(\bx)\}$ at each $\bx\in \calB_s$, i.e. $\Pi(\chi)$ is nonconvex on $\calB_s$.
It was shown by Gao and Ogden that these solutions are nonsmooth if  $\tau(\bx)$ changes its sign on $\calB_s$
 \cite{gao-ogden-qjmam}.

Analytical solutions for general 3-D finite deformation problem $(\calP)$ were first proposed  by Gao in 1998-1999 \cite{gao-ima98,gao-mecc99}.
It is proved recently \cite{gao-haj} that for St Venant-Kirchhoff material, the problem $(\calP)$ could have 24 critical  solutions at
each material point $\bx \in \calB$, but only one global minimizer.
The solution is unique if the external force is sufficiently large.

\begin{figure}[h]
  \includegraphics[width=14cm,height=5cm]{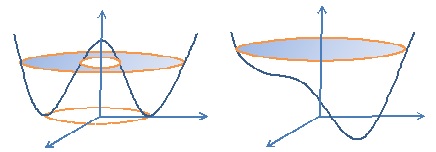}
 \caption{Graphs and level sets of $\PP(\bx)$ for  $ {\btau}= 0$ (left) and $ {\btau} \neq 0$ (right) }
 \label{1-dw} \vspace{.4cm}
 \end{figure}

 For a given function $\PP:\calF_a  \rightarrow \real $, its {\em level set } and  {\em sub-level set} of height $\alpha \in \real$
 are defined, respectively, as the following
  \eb
 {\cal L}_\alpha(\PP) = \{ \bF \in \calF_a\; | \;\; \PP(\bF) =  \alpha \}, \; \;    {\cal L}^\flat_\alpha(\PP) = \{ \bF \in \calF_a\; | \;\; \PP(\bF) \le \alpha \}, \; \;   \alpha \in \real.
  \ee
The geometrical explanation for ellipticity and Theorem \ref{thm-1} is illustrated by
 Fig. \ref{1-dw}, which shows that
 the nonconvex function
$\PP(\bgamma) = \half ( \half  |\bgamma|^2- 1)^2 - \bgamma^T \btau$ depends sensitively  on the
external force $\btau \in \real^2$.
If $|\btau |$ is big  enough,  $\PP(\bgamma) $
  has only one  minimizer and its
  level set  is an ellipse (Fig. \ref{1-dw} (b)).
  Otherwise,  $\PP(\bgamma)$ has multiple local minimizers and its level set is not an ellipse.
  For $\btau = 0$, it is  well-known  Mexican-hat in theoretical physics (Fig. \ref{1-dw} (a)). \vspace{1cm}

  Fig. \ref{1-dw}   shows    that although
  $\PP(\bgamma)$ has only one global minimizer for certain given ${\btau}$, the function is still nonconvex.
  Such a function is called quasiconvex in the context of global optimization.
   In order to distinguish this type of functions with  Morry's    quasiconvexity in nonconvex analysis,
   a generalized definition in a tensor space $\calF_a \subset \real^{m\times n}$ could be convenient.

  \begin{definition}[G-Quasiconvexity]
  A function $\PP:\calF_a \subset \real^{m\times n} \rightarrow \real $  is called G-quasiconvex
  if its domain $\calF_a$ is convex and
  \eb
  \PP( \theta \bF + (1- \theta) \bT) \le \max \{ \PP(\bF) , \PP(\bT) \} \;\;
  \forall \bF, \; \bT  \in \calF_a, \;\; \forall \theta \in [0,1].
  \ee
  It is called strictly G-quasiconvex if the inequality holds strictly.
\end{definition}
Moreover, we may need  a  definition of generalize   ellipticity   for nonconvex systems.
\begin{definition}[G-Ellipticity]
 For a given function $\PP:\calF_a \rightarrow \real$ and $\alp \in \real$, its
 level set ${\cal L}_\alp(\PP)$ is said to be a G-ellipse if it is a closed, simply connected set.
For a given $\bt $ such that $\bT \in \calT_a$, the $(BVP)$ is said to be $G$-elliptic if the total potential function
$\PP(\bF)$ is G-quasiconvex on $\calF_a$. $(BVP)$ is strongly G-elliptic if $\PP(\bF)$ is strictly G-quasiconvex.
\end{definition}

\begin{Lemma} For a given function $\PP:\calF_a \subset \real^{m\times n } \rightarrow \real $,
 \begin{eqnarray}
 \mbox{  $\PP(\bF)  $  is  G-quasiconvex  $\Leftrightarrow $  ${\cal L}^\flat_\alpha(\PP) $  is  convex
 $\Leftrightarrow $  ${\cal L}_\alpha(\PP) $  is a G-ellipse   $\forall \alpha \in \real$  .} \\
 \mbox{  $\PP(\bF)  $ is convex $\Rightarrow $ is rank-one convex $\Rightarrow $  G-quasiconvex  $\Rightarrow $  $(BVP)$ is G-elliptic} .
 \end{eqnarray}
 \end{Lemma}

This statement shows an important  fact  in nonconvex systems, i.e.
  the total number of solutions to a nonlinear equation depends not only on the stored energy, but also (mainly)  on  the external force field.
The nonlinear partial differential equation  in $(BVP)$
is elliptic only if it is  G-elliptic. $(BVP)$ has at most one  solution  if   $\PP(\bF)$
is strictly  G-quasiconvex on $\calF_a$.
\begin{rem}[Existence and Uniqueness] 
Suppose that the canonical function $\VV:\calE_a \rightarrow \real$ is convex, then $\nabla \VV^*(\zeta )> 0$ is a monotonic operator on
$\calE^*_a$. 
If for any  given   $\bt:\St \rightarrow \real^3$ such that $\bT \in \calT_a \neq \emptyset$ 
and $\tau^2 (\bx)= \tr( \bT^T \bT ) \neq  0 \;\; \forall \bx \in \calB$, then
the nonconvex variational problem $(\calP) $ has at least  one nontrivial solution a.e. in $\calB$.
It has a unique nontrivial  solution if  there exists a constant $\tau_c$ such that   $\tau^2 (\bx)= \tr( \bT^T \bT ) \ge \tau_c^2 \;\; \forall \bx \in \calB$.
\end{rem}


In global optimization,   the most simple
 quadratic integer programming problem
\[
(\calP)_i : \;\;\; \min \left\{ \Pi(\bx) = \half \bx^T \bQ \bx -  \bx^T {\bf t} \;\; | \;\; \bx = \{ x_i \}^n \in \{ 0, 1\}^n \subset \real^n \right\}
\]
could have up to $2^n$ local minimizers, which can't be solved   directly by traditional deterministic  methods  in polynomial time  due to the indefinite matrix $\bQ$ and the integer constraint.
Such a nonconvex discrete optimization problem is considered as NP-hard in computer science.
However, by using canonical transformation $\bxi =\Lambda(\bx)  = \{ x_i(x_i -1)\}\in \real^n$,    the canonical dual of this discrete problem
is  a concave maximization over a convex set in continuous space \cite{gao-bridge}. 
It was proved in \cite{gao-cace09} that there exists a positive vector $\btau = \{ \tau_i\}^n   > {\bf  0}  \in \real^n$,
if $\{ |t_i| \le \tau_i \}^n$, then   $\calS^+_a \neq \emptyset$ and $(\calP)_i$ is not NP-hard.
 The decision variable is simply $\{ x_i\}  = \{0  \mbox{ if }  t_i  <  - \tau_i , \; 1 \; \mbox{ if } t_i > \tau_i \} $ (Theorem 8, \cite{gao-cace09}).
 Thus, the canonical duality theory can be used to identify   NP-hard problems \cite{gao-bridge}.

\section{Applications in Anti-Plane Shear Deformation}
Now let us consider a special case that  the  homogeneous  elastic body
$\calB\subset \real^3$ is a cylinder
with generators parallel to the $\be_3$ axis and
with cross section a sufficiently nice region $\Oo \subset \real^2$ in the
$\be_1 \times \be_2$ plane.
The so-called anti-plane shear  deformation is defined by (see  \cite{horgan})
\begin{equation}\label{defor}
\bchi (\bx) = \left\{   x_1, \;\;   x_2, \;\;
   x_3  + \uu(\xx_1, \xx_2) \right\} : \Oo \subset \real^2 \rightarrow \real^3,
\end{equation}
where $ (\xx_1, \xx_2, \xx_3 )$ are cylindrical  coordinates in the
reference configuration  $\calB$ relative to a
cylindrical  basis $\{\be_i\},\; i = 1,2,3$.
 On $\Gu \subset \partial \Oo$, the homogenous boundary condition is given
$
 \uu(x_\alp)= 0 \;\; \forall x_\alp \in \Gu, \;\; \alp = 1,2 .$
On the remaining boundary $\Gt = \partial \Oo  \cap \Gu$, the cylinder is subjected to the  shear force
\[
{\bf t} (\bx) =  t(\bx)  \be_3 \;\; \forall \bx \in \Gt  ,
\]
where   $ t: \Gt \rightarrow \real $ is a prescribed function.
For this anti-plane shear deformation we have
\begin{equation}\label{1}
\mathbf{F}= \nabla \bchi =  \left(\begin{array}{ccc}
1  & 0 & 0 \\
0 &  1 & 0 \\
\uu_{,1} & \uu_{,2} & 1  \end{array} \right), \;\; \;\;
\mathbf{C}=\mathbf{F}^{\rm
T} \bF =
\left( \begin{array}{ccc}
  1 + \uu_{,1}^2  & \uu_{,1} \uu_{,2} &   \uu_{,1}\\
\uu_{,1} \uu_{,2}&  1 + \uu_{,2}^2 &   \uu_{,2} \\
  \uu_{,1} &   \uu_{,2} &  1
\end{array} \right) ,
\end{equation}
where  $\uu_{,\alpha}$ represents $\partial \uu/\partial x_\alpha$ for $\alpha= 1,2$.
By the notation $|\nabla \uu|^2 = \uu_{,1}^2 + \uu_{,2}^2$, we have
\begin{equation}
I_1(\bC)=  I_2 (\bC)= 3  + | \nabla \uu |^2, \;\;
  I_3 (\bC)  \equiv 1  , \la{Ii}
\end{equation}
Clearly, both $\bF$ and $ I_1(\bC)$ depend only on the shear strain $\bgamma = \nabla \uu = \{ \uu_{,\alp} \}$,
therefore, the strain energy can be equivalently written in the forms of
\eb
\WW(\bF(\bgamma)) = \VV(\xi(\bgamma)) = \hatW(\bgamma)
\ee
where $\hatW(\bgamma)$ is a real-valued function.

The fact   $\det \bF \equiv  1 $  shows that the anti-plane shear state  (\ref{defor})
  is an isochoric   deformation.
Therefore,
 the kinetically admissible displacement space $\calX_c$   can be simply replaced by a convex set
\eb
\calU_c = \{ \uu(\bx)\in \calW^{1,1}(\Oo; \real)| \;\; \uu(\bx)= 0 \;\; \forall \bx = \{ x_\alp\} \in \Gu \}. \label{eq-uc}
\ee
Thus, in terms of   $\xi =   \Lambda( \bgamma) = I_1 - 3 =  |\bgamma|^2 $
 and  $\WW(\bF(\bgamma)) = \barW(\Lambda(\bgamma)) $, for any given
 \[
 \btau \in  \calT_a =  \{ \btau \in C^1[\Oo; \real^2] | \;\; \nabla \cdot \btau = 0  \;\; \mbox{ in } \Oo,
\;\; \bn \cdot \btau = t \;\; \mbox{ on } \Gt \}
\]
Problem $(\calP)_{\bT}$ for the  anti-plane shear deformation  (\ref{defor})  has the following form
\eb
(\calP)_s: \;\; \min \left\{ \Pi(\uu) = \int_\Oo \PP(\nabla\uu) \dO \;\;| \;\; \uu \in \calU_c \right\}, \;\;
\PP(\bgamma) = \barW(\Lambda(\bgamma)) - \bgamma^T \btau
\ee

Under certain regularity conditions, the associated  mixed boundary value problem is
\eb\label{eq-eqs}
(BVP)_s: \;\;    \left\{ \begin{array}{l}
\nabla \cdot \left( 2  \zeta   \nabla  \uu  \right)  = 0 \;\; \mbox{ in } \Oo, \\
\bn \cdot \left( 2  \zeta   \nabla  \uu  \right)= t \;\; \mbox{ on } \Gt, \;\; \uu = 0 \;\; \mbox{ on } \Gu
\end{array} \right.
\ee
where  $ \bn =  \{ n_\alpha\} \in \real^2$ is a unit vector  norm to $\partial \Oo$, and
$ \zeta  = \nabla  \barW(\xi), \;\;  \xi =  |\nabla \uu|^2$.

If $\Gu = \partial \Oo$, then $(BVP)_s$ is a Dirichlet boundary value problem, which has only trivial solution due to zero input. For Neumann boundary value problem $\Gt = \partial \Oo$, the external force field must be such that
\[
\int_\Gt t(\bx) \dG  = 0
\]
for overall force equilibrium. In this case, if $\barbchi$ is a solution to  $(BVP)_s$, then $\bchi = \barbchi + {\bf c}$ is also
a solution for any vector ${\bf c} \in \real^3$ since the cylinder is not fixed. Therefore, the mixed boundary value problem $(BVP)_s$ is necessary for anti-plane shear deformation to have a unique solution.

By the fact that the only unknown $\uu$ is a scalar-valued function,
 anti-plane shear deformations are one of the simplest classes of deformations that solids can undergo \cite{horgan}.
Indeed, if $\barW(\xi)$ is a canonical function on
 $\calE_a = \{ \xi \in L^p(\Oo) | \;\; \xi (\bx) \ge 0 \;\; \forall \bx \in \Oo\}$  and for any given
$ \btau \in \calT_a$
 such that $\tau= |\btau|$,
 the canonical dual problem has a very simple form
  \eb
 (\calP^d)_s: \;\;\;\; \sta \left\{  \Pi^d(\vsig) =   \int_{\Oo}  \left[
   - \barW^*(\vsig) - \frac{1}{4} \vsig^{-1}  \tau^2  \right] \dO \; | \;\; \vsig \in \calS_a \right\} .
  \ee
Since  $\Lambda(\uu) = |\nabla \uu|^2  $, the canonical dual algebraic equation (\ref{cda}) for this problem is
\eb
4 \zeta^2  \nabla \VV^*(\zeta)  = \tau^2(\bx),\;\;\;\; \forall \bx \in \Oo . \label{cda1}
\ee
\begin{Corollary}\label{coro}
For any given  non-trivial shear force $ t (\bx)   \neq 0 $ on $ \Gt$ such that $\btau \in \calT_a \neq \emptyset$, the canonical dual problem $(\calP^d)_s$  has at least one non-trivial solution   $\vsig_k$. If $\nabla \times \vsig_k^{-1} \btau = 0$, the scale-valued function
   \eb\label{eq-asolu}
   \uu_k (\bx) = \half \int_{\bx_0}^\bx  \vsig_k^{-1}  \btau \cdot \mbox{d}\bx
   \ee
   along any path from $\bx_0 \in \Gu$ to $\bx \in \Oo$ is a critical point of $\Pi(\uu)$ and
   $\Pi(\uu_k) = \Pi^d( \vsig_k) $.

   If $\zeta_k \in \calS^+_a$, then $\uu_k$ is a global minimizer of $(\calP)_s$.

    If $\zeta_k \in \calS^-_a$ and $\bG(\zeta_k) \succ 0$, then $\uu_k$ is a local  minimizer of $(\calP)_s$.

      If $\zeta_k \in \calS^-_a$ and $\bG(\zeta_k) \prec 0$, then $\uu_k$ is a local  maximizer of $(\calP)_s$.

 \end{Corollary}

\begin{figure*}[htbp]
  \centering
\subfigure[]{
    \label{ } 
    \includegraphics[width=2.4in]{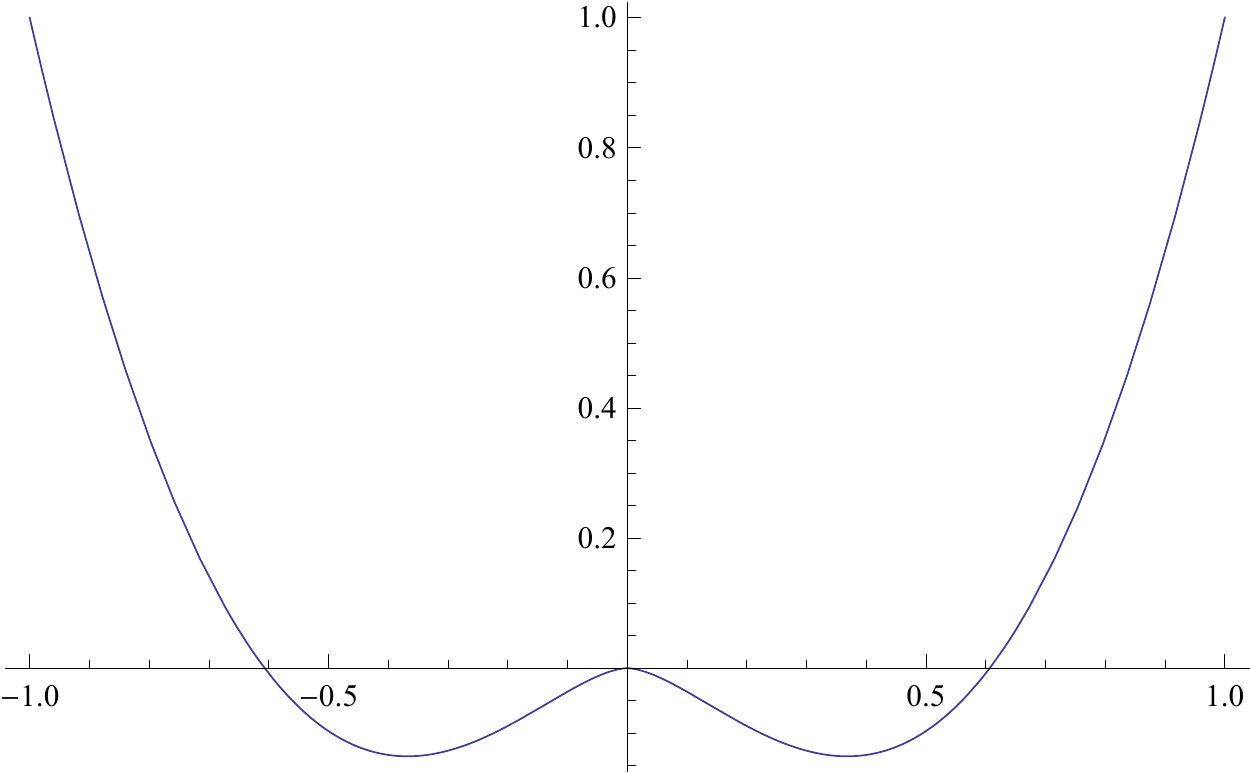}}
    \hspace{0.1in}
\subfigure[ ]{
    \label{ } 
    \includegraphics[width=2.4in]{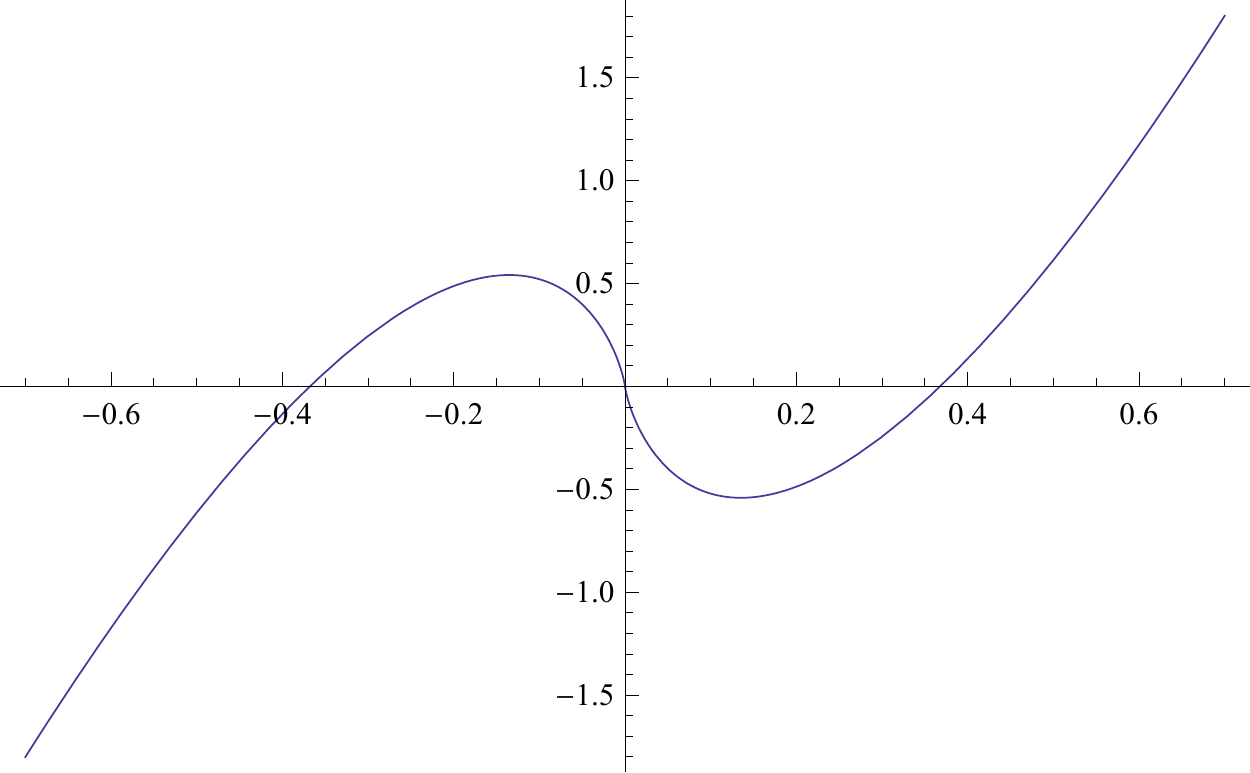}}
\caption{Graphs of $\hat{W}(\bgamma)$ (a) and its derivative (b) ($\cc_1= \cc_2=1$)}
  \label{dw1d} 
\end{figure*}

\noindent {\bf Example.} Applications of the canonical duality theory  to general anti-plane shear problems have been demonstrated
for solving convex exponential and nonconvex polynomial stored energies recently in \cite{gao-cmt15}.
In this paper, the following generalized neo-Hookean model is considered  
\eb
\VV(\xi) = \cc_1 (I_1 - 3) + \cc_2 (I_1-3)  \log(I_1 -3)
\ee
where  $\cc_1, \; \cc_2$ are positive material constants.
Clearly, $\VV(\xi)$ is convex  in $\xi= I_1 -3 $, but
\[
\hatW(\bgamma) = \VV(I_1(\bgamma)) = \cc_1 |\bgamma |^2 + \cc_2 |\bgamma|^2 \log |\bgamma|^2
\]
 is a
double-well function of the shear strain $\bgamma = \nabla \uu$ (see Fig. \ref{dw1d}).

It is easy to check
\[
\vsig = \nabla \VV(\xi) = \cc_1 + \cc_2 (\log \xi + 1) : \calE_a \rightarrow \calE_a^*
=  L^q(\Oo)
\]
is one-to-one and onto, where $q$ is a dual number of $p\ge 1$, i.e. $1/p + 1/q =1$. The complementary energy can be obtained easily
\[
\VV^*(\zeta) = \sta  \{ \xi \zeta - \VV(\xi) | \;\; \xi \in \calE_a \} = \cc_2 \exp [\cc_2^{-1}(\zeta - \cc_1) - 1]
\]
In this case, the canonical dual algebraic equation is
\eb \label{eq-dab3}
\zeta^2 \exp\left[ \frac{\zeta - \cc_1}{\cc_2} - 1 \right] = \tau^2(\bx) \;\; \forall \bx \in \Oo.
\ee
  Let $h^2(\zeta)=   \zeta^2     \exp [(\zeta - \cc_1)/\cc_2 - 1 ]  $
  be the left hand side function in
 the canonical  dual algebraic equation (\ref{eq-dab3}). By solving $h'(\zeta_c) = 0$ we known that  at  $\zeta_c= - 2 \cc_2$, $h(\zeta)$ has a local maximum
 \[
 h_{\max}(\zeta_c)  = \eta =   2 \cc_2 \sqrt{\exp[-3 - \cc_1/\cc_2]} .
 \]

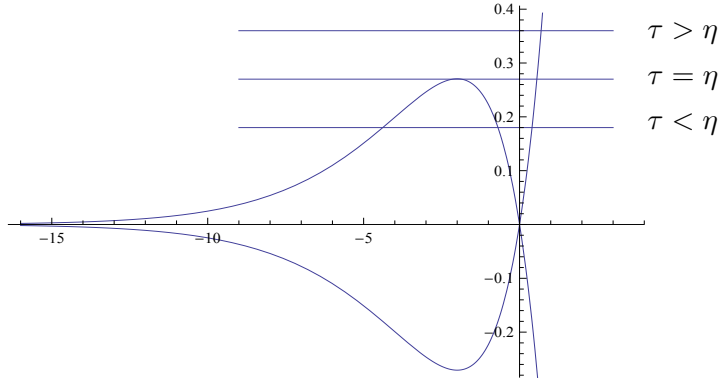
\begin{figure}[h]
\setlength{\unitlength}{.3cm}
\begin{picture}(-1,12)
\put(14,15){\line(1,0){15}}
\put(33,15){$\tau  >  \eta $}
\put(14,13.2){\line(1,0){15}}
\put(33,13){$\tau   =\eta $}
\put(14,11){\line(1,0){15}}
\put(33,11){$\tau  < \eta $}
 \end{picture}
  \hspace{.6in} \includegraphics[width=3.3in,height=5cm]{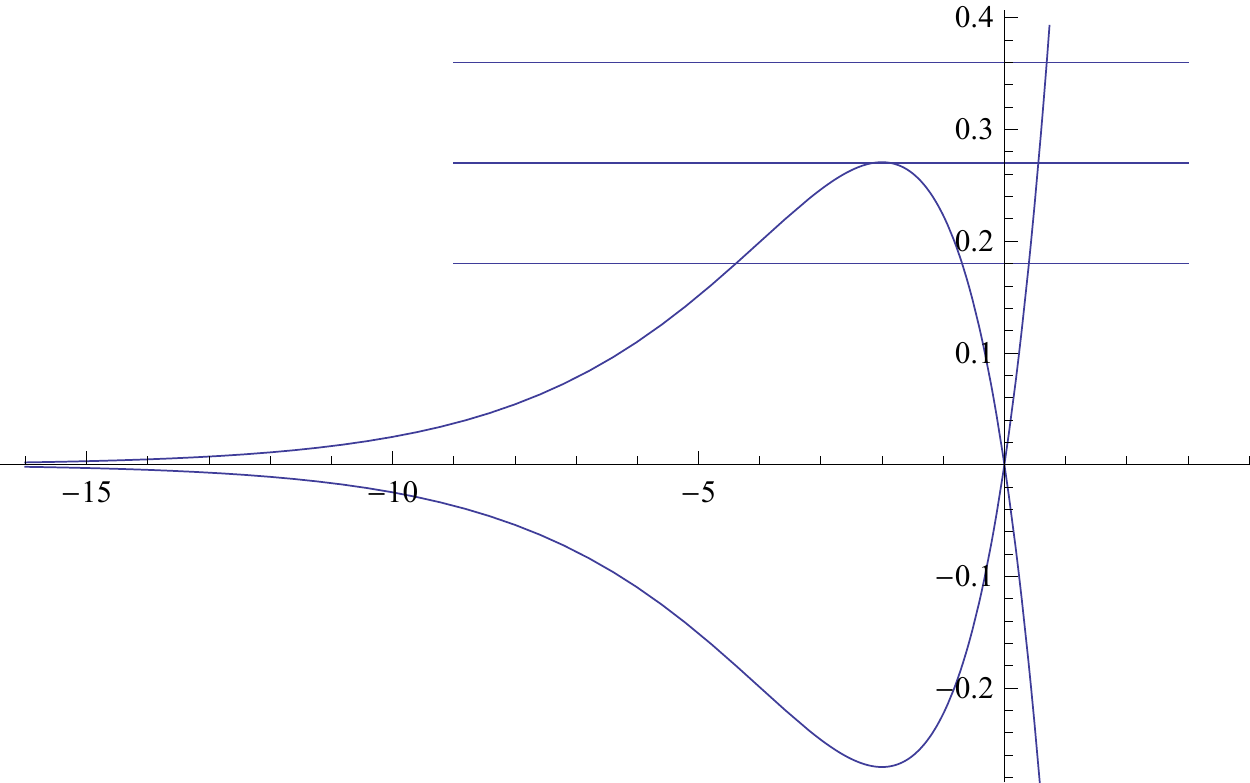}
\caption{Dual algebraic curve   $ h(\zeta)  \;( \cc_1 =   \cc_2 = 1) $}\label{fig-admifig3}
\end{figure}

\begin{figure*}[htbp]
  \centering
\subfigure[$\tau  > \eta$]{
    \label{fig5:subfig:a} 
    \includegraphics[width=2.4in]{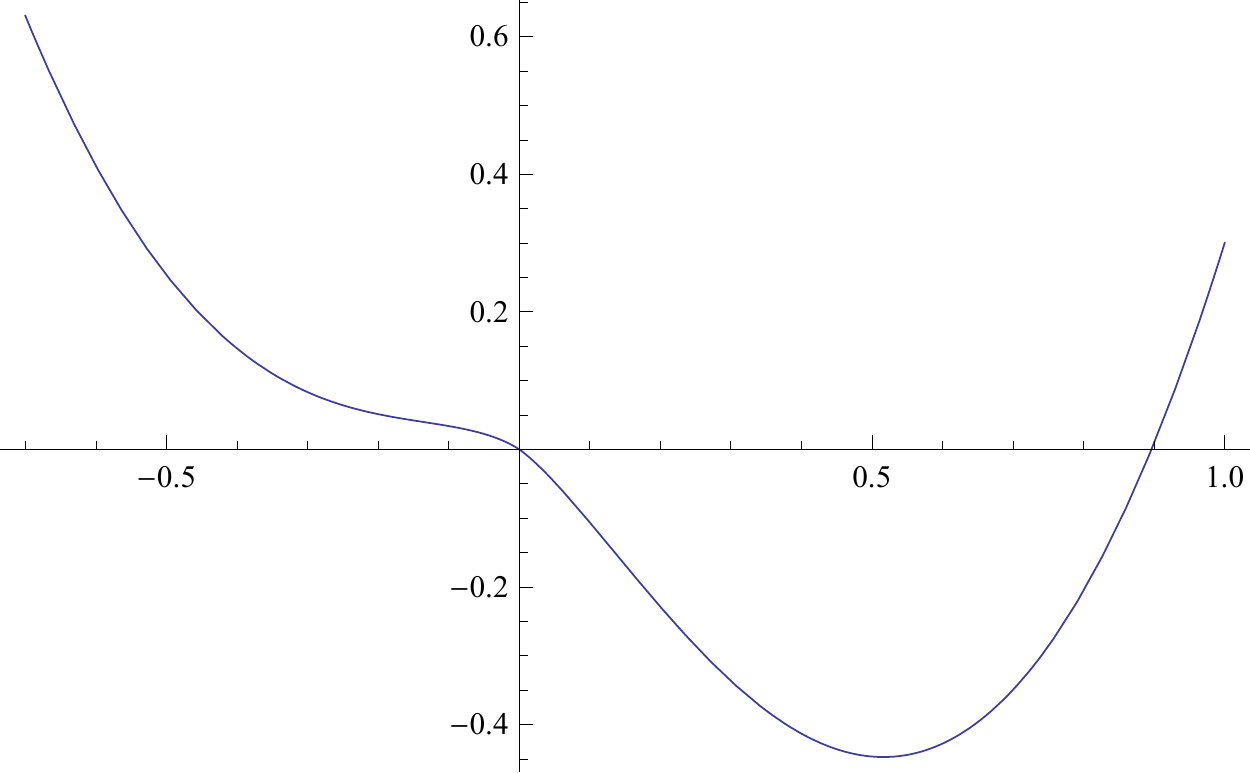}}
    \hspace{0.1in}
\subfigure[ $\tau  = \eta$  ]{
    \label{fig5:subfig:b} 
    \includegraphics[width=2.4in]{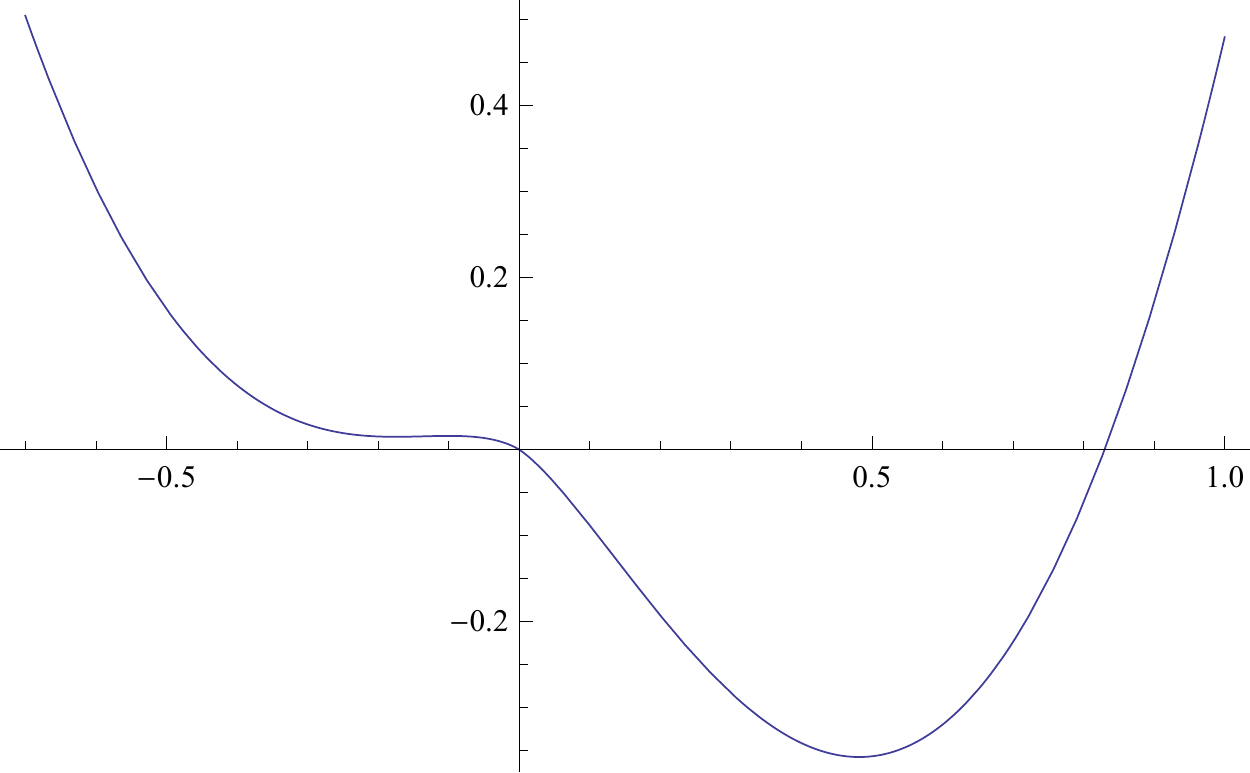}}
\caption{Graphs of G-quasiconvex  $ G(\gamma)$  ($\cc_1=  \cc_2=1$)}  \label{fig5}
\end{figure*}

\begin{figure*}[htbp]
  \centering
\subfigure[Graph  of $G(\gamma)$]{
    \label{fig6:subfig:a} 
    \includegraphics[width=2.8in]{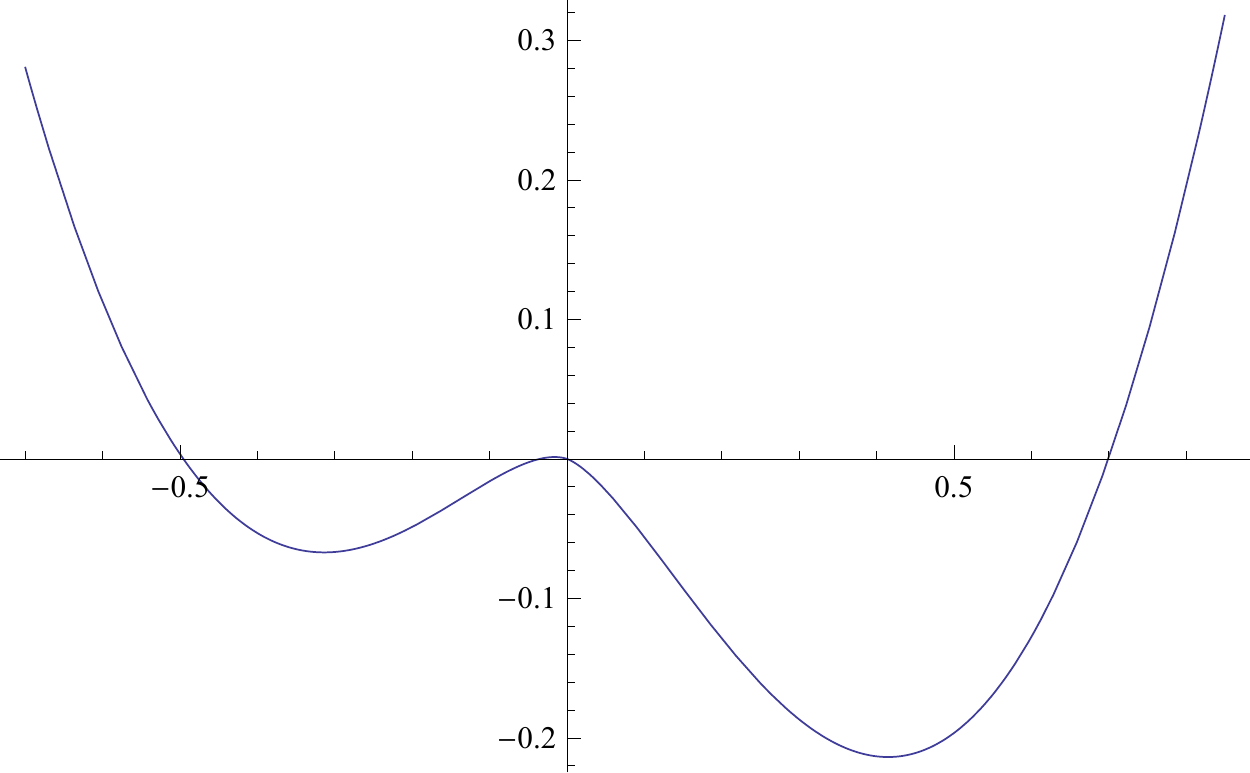}}
    \hspace{0.1in}
\subfigure[Graph  of $G^d(\zeta)$]{
    \label{fig6:subfig:b} 
    \includegraphics[width=2.8in]{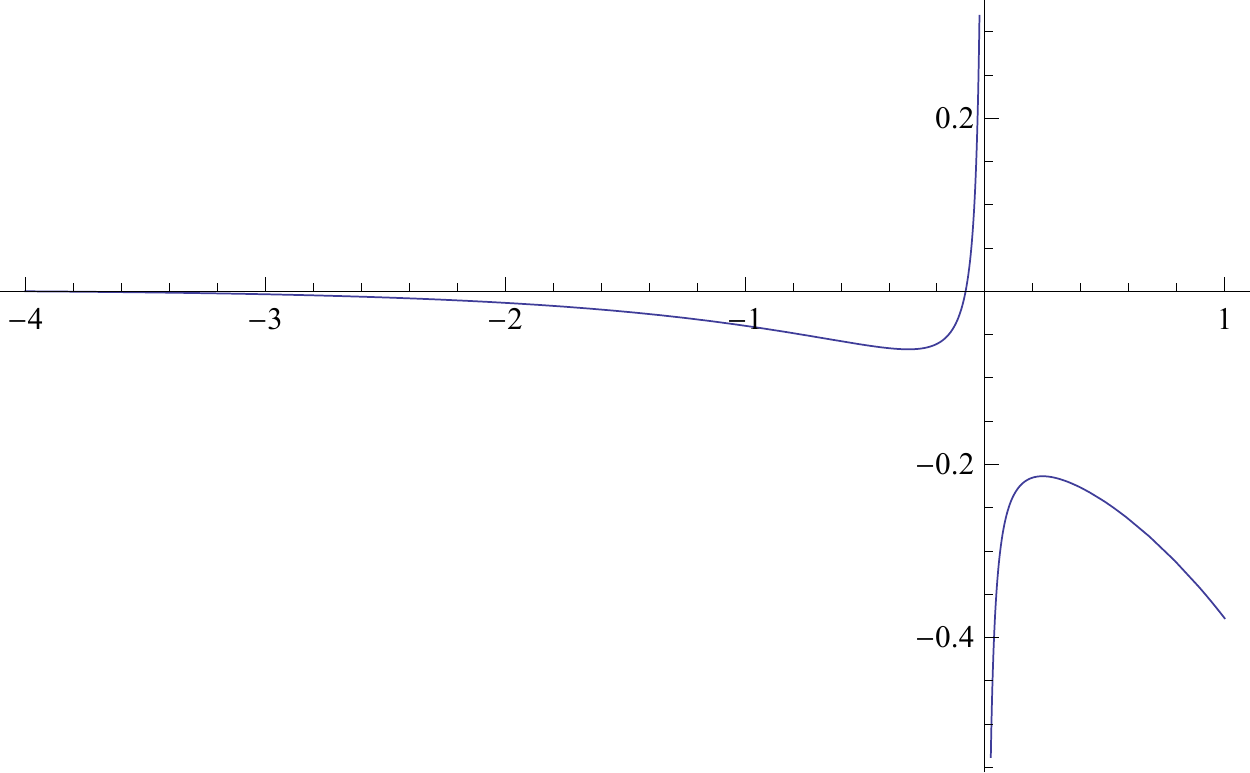}}
\caption{Graphs of $ G(\gamma)$ and $ G^d(\zeta)$ for  $\tau  <  \eta$    ($\cc_1= \cc_2 =1 $ )}  \label{fig6} 
\end{figure*}

  From the graphs of the canonical  dual algebraic curve
  $ h(\zeta) $  
   given in Fig. \ref{fig-admifig3}   we can see that
   the canonical dual algebraic equation (\ref{eq-dab3}) may have at most
   three real  solutions in the order of $\zeta_1 \ge 0 \ge \zeta_2 \ge \zeta_3$
  depending on  $\tau = |\btau(\bx)|  , \;\; \bx \in \Oo$
  (see Fig. \ref{fig-admifig3}b).
The equation  (\ref{eq-dab3}) has a unique solution if $\tau > \eta$. In this case, the total strain grand
$ G(\gamma)$ is strictly G-quasiconvex (see Fig. \ref{fig5}).
Fig \ref{fig6} shows the graphs of $G(\gamma)$ and its canonical dual $G^d(\zeta)$ for $\tau  < \eta$.
In this case, the function $G(\gamma)$ is nonconvex and has three critical points. The triality theory holds for
$G(\gamma)$ and its canoncial dual $G^d(\zeta)$
\[
 G(\gamma_1) = \min_{\gamma \ge 0} G(\gamma) = \max_{\zeta > 0}  G^d(\zeta) = G^d(\zeta_1).
 \]
 \[
 G(\gamma_2) = \min_{\gamma \in \calG_o} G(\gamma) = \min_{\zeta  > -2\cc_2}  G^d(\zeta) = G^d(\zeta_2).
\]
\[
G(\gamma_3) = \max_{\gamma \in \calG_o } G(\gamma) = \max_{\zeta < - 2 \cc_2}  G^d(\zeta) = G^d(\zeta_3),
\]
where $\calG_o $ is a neighborhood    of $\gamma_i $ $(i=1,2)$. \vspace{-.0cm}

\section{Conclusions}
In summary,  the following conclusions can be   obtained. 
\begin{verse}
1.  The pure complementary energy principle and canonical duality-triality  theory developed in \cite{gao-dual00} are useful
 for
 solving general nonlinear boundary value problems in nonlinear elasticity. \\

 2. Both convexity of the total potential and  ellipticity condition of the associated fully nonlinear boundary value problem
 depend  not only on the stored energy function,
 but also sensitively on the external force field.\\

 3.   The Legendre-Hadamard  condition is only a necessary ellipticity condition for convex systems.
 The   triality theory provides a sufficient condition to identify  both global and local extremum solutions for nonconvex problems.\\


  \end{verse}

These  conclusions are naturally included in the canonical duality-triality theory developed by the author and his co-workers during the last 25 years \cite{gao-dual00}.
Extensive applications have been given in multidisciplinary fields of biology, chaotic dynamics, computational mechanics, information theory, phase transitions, post-buckling,  operations research, industrial and systems engineering, etc.  (see  recent review article \cite{gao-bridge}).

\subsection*{Acknowledgements}
Insightful discussions with  Professor David Steigmann from  UC-Berkeley is  sincerely acknowledged.
The research   was supported by US Air Force Office of Scientific Research (AFOSR FA9550-10-1-0487).

\end{document}